\documentclass{amsart}

\usepackage{amsfonts}

\usepackage{amsmath}
\usepackage{amssymb}
\usepackage{amsthm}
\usepackage{latexsym}

\theoremstyle{plain}
\newtheorem{theorem}{Theorem}

\theoremstyle{definition}

\theoremstyle{remark}

\newcommand{\set}[1]{\left\{#1\right\}}

\begin{document}

\title{Mapping the Discrete Logarithm} 

\author{Daniel R. Cloutier}
\address{Cincinnati, OH, USA}
\email{Daniel.R.Cloutier@alumni.rose-hulman.edu}

\author{Joshua Holden}
\address{Rose-Hulman
Institute of Technology, Terre Haute, IN 47803, USA}
\email{Joshua.Holden@rose-hulman.edu}

\begin{abstract}
The discrete logarithm is a problem that surfaces frequently in the
field of cryptography as a result of using the transformation
$g^a\mod n$.  This paper focuses on a prime modulus, $p$, for which
it is shown that the basic structure of the functional graph is
largely dependent on an interaction between $g$ and $p-1$.  In fact,
there are precisely as many different functional graph structures as
there are divisors of $p-1$.  This paper extracts two of these
structures, permutations and binary functional graphs.  Estimates
exist for the shape of a random permutation, but similar estimates
must be created for the binary functional graphs.  Experimental data
suggests that both the permutations and binary functional graphs
correspond well to the theoretical data which provides motivation to
extend this to larger divisors of $p-1$ and study the impact this
forced structure has on the many cryptographic algorithms that rely
on the discrete logarithm for their security. This is especially
applicable to those algorithms that require a ``safe" prime
($p=2q+1$, where $q$ is prime) modulus since all non-trivial
functional graphs generated using a safe prime modulus can be
analyzed by the framework presented here.
\end{abstract}

\maketitle

\section{Introduction}\label{S:Intro}
Just a few decades ago, cryptography was considered a domain
exclusive to national governments and militaries.  However, the
computer explosion has changed that. Every day, millions of people
trust that their privacy will be protected as they make online
purchases or communicate privately with a friend.  Many of the
cryptographic algorithms they will use are built upon a common
transformation, namely
\begin{equation}\label{EQ:General}
g^x\equiv y\mod{n}.
\end{equation}
For instance, Diffie-Hellman key exchange, RSA and the Blum-Micali
pseudorandom bit generator all use (\ref{EQ:General}).  This paper
will examine some of the properties exhibited by this sort of
transformation and provide theoretical and experimental data
describing how the interaction between $g$ and the modulus impacts
the behavior of this function. 

\section{Terminology and Background}\label{S:TermAndBack}

In this paper, we will restrict the values of $n$ to primes and
examine mappings 
$$f : S = \set{1,2,...,p-1} \to S$$
of the form $x \mapsto g^{x} \bmod{p}$, where $p$ is a prime modulus.
In some instances, it will prove to be useful to interpret the
mappings as functional graphs.  A functional graph is a directed graph
such that each vertex must have exactly one edge directed out from it.
The relationship between the mappings which interest us and functional
graphs is straightforward.  Each element in $S$ can be interpreted as
a vertex.  The edges are defined such that an edge $\langle
a,b\rangle$ is in the graph if and only if $f(a) = b$.

There are a number of statistics of interest derived from
functional graphs. Following the convention of~\cite{RandMap},
which treats random mappings in detail, let $f:S\rightarrow
S$ be the transition function so that the edges in the functional
graph can be expressed as the ordered pair $\langle
x,f(x)\rangle$ for $x,f(x)\in S$. By applying the
pigeonhole principle and noting that the cardinality of $S$ is
$p-1$ we can say that by starting at any random point $u_0$ and
following the sequence $u_1 = f(u_0)$, $u_2 = f(u_1)$,
..., there must be a $u_i = u_j$ after at most $p$ iterations.
Suppose $u_i$ occurs before $u_j$ in the sequence of nodes. In
this case, the tail length is the number of iterations from $u_0$
to $u_i$. The cycle length is the number of iterations from $u_i$
to $u_j$. In more natural graphical terms, the cycle length is the
number of edges (or equivalently nodes) involved in the directed
path from $u_i$ to itself.  The tail length is the number of edges
from $u_0$ to $u_i$. Additionally, a terminal node is one with no
pre-image, or more formally, $x$ is a terminal node if
$f^{-1}(x) = \emptyset$. A node is an image node if it is
not a terminal node. Since each node has an out-degree of exactly
one, each cycle with the trees grafted onto its nodes will form a
connected component.

The value of $g$ plays a major role in determining the basic
structure of the graph.  In fact, as Theorem~\ref{T:m-ary}
formalizes, the interaction among $g$ and $p-1$ will effectively
fix the in-degrees of the nodes in the graph.  First, though,
define an $m$-ary functional graph to be a graph where each node
has in-degree of exactly zero or $m$.  The proof of the following 
theorem is then straightforward.

\begin{theorem}\label{T:m-ary}
Let $p$ be fixed and let $m$ be any positive integer that divides
$p-1$.  Then as $g$ ranges over all integers, there are
$\phi(\frac{p-1}{m})$ different functional graphs which are $m$-ary
produced by maps of the form $f: x \mapsto g^{x} 
\bmod{p}$.\footnote{Throughout this paper, $\phi$ denotes the Euler 
phi funcion.}
Furthermore, if $r$ is any primitive root modulo $p$, and $g\equiv
r^a\mod p$, then the values of $g$ that produce an $m$-ary graph are
precisely those for which $\gcd(a,p-1) = m$.
\end{theorem}

Theorem~\ref{T:m-ary} gives a strong indication that the graphs
generated by~\eqref{EQ:General} have to be considered separately
for different values of $m$.  

It should be noted, though, that there are
some values of $m$ which lead to completely predicable graphs. For
instance, there is one $(p-1)$-ary graph that corresponds to
$g\equiv 1\mod p$.  There is also one $(\frac{p-1}{2})$-ary graph
that corresponds to $g\equiv -1\mod p$.  In general, however, an
$m$-ary graph is not trivially predictable.  This paper will
restrict its focus to unary functional graphs (which will be
referred to as permutations since they simply permute the numbers
1, ..., $p-1$) and binary functional graphs. The values of $g$
which produce a permutation are precisely those which are
primitive roots modulo $p$.

In cryptography, it is common to look for primes where $p-1$ has at
least one large prime factor.  For instance, the pseudorandom bit generator described by
Gennaro in~\cite{PRBG} and mentioned in Section~\ref{S:Intro}
requires the modulus to be of the form $p = 2q+1$ where $q$ is also
prime.  A prime of this form is known as a safe prime ($q$ is also
known as a Sophie Germain prime). These primes are of interest here
not only because of their extensive use in cryptography, but also
because $p-1$ has only four divisors, namely 1, 2, $q$ and $2q$. It
can be quickly verified that there is only one $q$-ary ($g\equiv
-1\mod{p}$) and one $2q$-ary ($g\equiv 1\mod{p}$) graph generated.
More importantly, there are $\phi(q)$ permutations and $\phi(q)$
binary functional graphs which represent the remaining values of $g$
(since $\phi(q)$ is $q-1$). Thus, not only do safe primes provide
large numbers of permutations and binary functional graphs, but
every graph generated by a safe prime is either trivial (the graphs
where $g$ is either 1 or -1) or fits into the theoretical framework
presented in Section~\ref{S:Theory}.

\section{Theoretical Results}\label{S:Theory}

In Theorem~\ref{T:m-ary}, it is shown that the in-degree of each node
is dependent on the value of both $g$ and $p$.  This is clearly
imposing a structure on any functional graphs generated
using~\eqref{EQ:General}.  It seems reasonable, though, that a large
collection of functional graphs generated by using~\eqref{EQ:General}
as the transition function would tend toward exhibiting behavior
similar to that of a collection of random functional graphs.  At a
minimum, a factorization for $p-1$ with many divisors would
certainly seem to hide the structure imposed by Theorem~\ref{T:m-ary}
since the many divisors of $p-1$ would each contribute some graphs.
Section~\ref{SS:CombinedObserved} will give evidence that this is not
the case.  However, the methods used to obtain the theoretical bounds
for the random functional graphs can be extended to analyze $m$-ary
graphs for specific $m$.

While most of the parameters that are of interest depend on the
exact graph generated, the number of image nodes can be computed
directly from the values of $g$ and $p$.  The proof is again 
straightforward.

\begin{theorem}\label{T:image}
The number of image nodes in any $m$-ary graph is $\frac{p-1}{m}$.
\end{theorem}

Theorem~\ref{T:image} helps to quantify the repercussions of
Theorem~\ref{T:m-ary} and the restrictions on in-degree in $m$-ary
graphs.  The number of image nodes is a direct function of $m$
which can greatly limit the shapes each graph can take on.  None
of the other parameters appear to have a generalization as
convenient as the image nodes and will be treated as specific
parameters in permutations and binary functional graphs.

\subsection{Random Functional Graphs}\label{SS:Unconstrained}

Flajolet and Odlyzko do a thorough analysis of functional graphs
in~\cite{RandMap}. While none of these results are original,
Flajolet and Odlyzko demonstrate that all of these parameters can
be estimated through a singularity analysis of generating
functions.  This appears to be the first method that can be
applied to all of these parameters. Their methods can then be
adapted for any fixed value of $m$ to estimate the parameters of
interest for an $m$-ary graph. Specifically, the methods will be
used to confirm some permutation results and to develop all of the
binary functional graph results. The results from~\cite{RandMap}
are summarized below in Theorem~\ref{T:FnGraphs}.
\begin{theorem}\label{T:FnGraphs}
The asymptotic values for the parameters of interest in a random
functional graph of size $n$ are:
\begin{align}
\textit{Number of components }\quad &
\frac{\ln{(2n)}+\gamma}{2}\tag{i}\label{T:FnGraphComp}\\
\textit{Number of cyclic nodes }\quad &
\sqrt{\pi n/2}-\frac{1}{3}\tag{ii}\label{T:FnGraphCycNodes}\\
\textit{Number of tail nodes }\quad &
n-\sqrt{\pi n/2}+\frac{1}{3}\tag{iii}\label{T:FnGraphTailNodes}\\
\textit{Number of terminal nodes }\quad &
e^{-1}n\tag{iv}\label{T:FnGraphTerminal}\\
\textit{Number of image nodes }\quad &
(1-e^{-1})n\tag{v}\label{T:FnGraphImage}\\
\textit{Average cycle length }\quad &
\sqrt{\pi n/8}\tag{vi}\label{T:FnGraphAvgCycle}\\
\textit{Average tail length }\quad &
\sqrt{\pi n/8}\tag{vii}\label{T:FnGraphAvgTail}\\
\textit{Maximum cycle length }\quad &
\sqrt{\frac{\pi n}{2}}\int_0^\infty\left[1-\exp{\left(-\int_v^\infty e^{-u}\frac{du}{u}\right)}\right]\,dv\notag\\
&\approx 0.78248\sqrt{n}\tag{viii}\label{T:FnGraphMaxCycle}\\
\textit{Maximum tail length }\quad &
\sqrt{{2 \pi n}}\ln 2 \approx 1.73746\sqrt{n}\tag{ix}\label{T:FnGraphMaxTail}\\
\end{align}
\end{theorem}

In part~(\ref{T:FnGraphComp}), $\gamma$ refers to the Euler constant
which is approximately 0.57721566.  The second order terms for
parts~(\ref{T:FnGraphComp}),~(\ref{T:FnGraphCycNodes}),
and~(\ref{T:FnGraphTailNodes}) were not given in~\cite{RandMap}, but
can be computed with a careful singularity analysis using precisely
the same methods used there.

\subsection{Permutations}\label{SS:PrimRoot}
Predicting the behavior of the permutations is, in many ways, much
easier than other $m$-ary graphs.  The most important reason for
this is that there are no terminal nodes or tail nodes.  This
follows quickly from the definition of a permutation as a unary
functional graph and the fact that the sum of the in-degrees must
be the same as the sum of the out-degrees.  Each node has an
out-degree of exactly one, and if any node were to have an
in-degree of zero, then, by the pigeon-hole principle, at least
one node must have an in-degree of more than one.  This is not
allowed so each node must have in-degree of exactly one.
Furthermore, since every tail must contain at least one terminal
node, this also implies that every node is cyclic.  The parameters
that can then be determined from the definition of a permutation
are given below. 
\begin{align}
&&&\text{Number of cyclic nodes} &&n&&\notag\\
&&&\text{Number of tail nodes} &&0&&\notag\\
&&&\text{Number of terminal nodes} &&0&&\notag\\
&&&\text{Number of image nodes} &&n&&\notag\\
&&&\text{Average tail length} &&0&&\notag
\end{align}
There are three non-trivial parameters of interest.  They are
expressed in Theorem~\ref{T:Permutation}.
\begin{theorem}\label{T:Permutation}
The asymptotic values for the number of components, the average
cycle length as seen from a random node and the maximum cycle length
in a random permutation of size $n$ have the following values:
\begin{align}
\textit{Number of components }\quad&\sum_{i=1}^{n} \frac{1}{i}\tag{i}\label{T:PermComp}\\
\textit{Average cycle length }\quad&\frac{n+1}{2}\tag{ii}\label{T:PermAvgCyc}\\
\textit{Maximum cycle length }\quad&
n\int_0^\infty\left[1-\exp{\left(-\int_v^\infty e^{-u}\frac{du}{u}\right)}\right]\,dv\notag\\
&\approx0.62432965n\tag{iii}\label{T:PermMaxCyc}
\end{align}
\end{theorem}

Parts~(\ref{T:PermComp}) and~Part~(\ref{T:PermAvgCyc}) are fairly 
well known.  Part~(\ref{T:PermMaxCyc}) seems to have first been solved by Shepp
and Lloyd in 1966~\cite{SheppLloyd66}.  An alternative solution
and proof more similar to the methods used here is offered by
Flajolet and Odlyzko in~\cite{SingAnalysis}.

\subsection{Binary Functional Graphs}\label{SS:SafePrime}
While estimates for the parameters investigated here exist in
literature for the random functional graphs and permutations, it
does not appear similar estimates exist for binary functional
graphs.  However, the methods in~\cite{RandMap}  can be extended to
develop these estimates.  Imitating the methods of~\cite{RandMap},
we first need to convert our ideas of a binary functional graph
into corresponding generating functions. The machinery is fairly
straightforward once we define the following as in~\cite{RandMap}:
\begin{center}
\begin{tabular}{ll}
BinFunGraph & = set(Components)\\
Component & = cycle(Node*BinaryTree)\\
BinaryTree & = Node + Node*set(BinaryTree, cardinality = 2)\\
Node & = Atomic Unit
\end{tabular}
\end{center}
This implies that a binary functional graph is a set of
components. Each component is a cycle of nodes with each node
having an attached binary tree to bring its in-degree to two.  A
binary tree is either a node (terminal node) or a node with two
binary trees attached. Finally, a node is simply an atomic unit.
A moment's reflection should indicate that this natural
specification does, in fact, specify a binary functional graph.
Imitating the transformations in~\cite[Section~2.1]{RandMap}, the
generating functions of interest are
\begin{align}
f(z) = e^{c(z)} = \frac{1}{1-zb(z)}\label{EQ:F}\\
c(z) = \ln{\frac{1}{1-zb(z)}}\label{EQ:C}\\
b(z) = z + \frac{1}{2}zb^2(z)\label{EQ:B}
\end{align}
Here $f$ generates the number of binary functional graphs, $c$
generates the number of components, and $b$ generates the number of
binary trees of a given size.  Solving the quadratic formula
for~(\ref{EQ:B}), we can produce the following formulas for $f$ and $c$
which simplify some of the cases:
\begin{align}
f(z) = \frac{1}{\sqrt{1-2z^2}}\label{EQ:F2}\\
c(z) = \ln{\frac{1}{\sqrt{1-2z^2}}}\label{EQ:C2}
\end{align}
In order to compute asymptotic forms of any of the statistics of
interest, we must first compute an asymptotic form for $f$ to
normalize results.  The following derivations give only a highlight of
the methods used by Flajolet and Odlyzko.  The interested reader is
encouraged to see~\cite{SingAnalysis,RandMap} for detailed proofs.

From equation~(\ref{EQ:F2}) it is clear that there is a singularity at
$z=1/\sqrt{2}$.  Performing a singularity analysis\footnote{The
analyses in this paper have been performed using the computer algebra
program Maple and the packages created as part of the Algorithms
Project at INRIA, Rocquencourt, France.  The packages can be found
online at http://pauillac.inria.fr/algo/libraries/software.html.} as
in~\cite[Section~2]{RandMap}, the asymptotic form for $f$ falls out
quickly as
\begin{equation}\label{EQ:NumGraphs}
f(z) \sim \frac{2^{n/2}}{\sqrt{\pi n/2}}.
\end{equation}
In at least one case, there are some important second-order interactions
between the error terms of the number of graphs and the appropriate
statistic. In these cases, a more exact form of~(\ref{EQ:NumGraphs})
must be used.  Expanding one more term in the expansion of $f$
gives
\begin{equation}\label{EQ:NumGraphsExact}
f(z) \sim \frac{2^{n/2}}{\sqrt{\pi n/2}} -
\frac{2^{n/2}}{4n\sqrt{\pi n/2}} =
\frac{2^{n/2}(4n-1)}{4n\sqrt{\pi n/2}}
\end{equation}
In most cases, using this more precise expansion of $f$ is not
necessary and does not change the results. Therefore, in all but
the necessary cases,~(\ref{EQ:NumGraphs}) will be used.

We begin by deriving the results for the most simple parameters.
\begin{theorem}\label{T:BinFunDirect}
The asymptotic forms for the number of components, number of
cyclic nodes, number of tail nodes, number of terminal nodes and
number of image nodes in a random binary functional graph of size
$n$, as $n\rightarrow\infty$ are
\begin{align}
\textit{Number of components }\quad &\frac{\ln{(2n)}+\gamma}{2}
\tag{i}\label{T:BinFnGraphComp}\\
\textit{Number of cyclic nodes }\quad & \sqrt{\pi n/2}-1
\tag{ii}\label{T:BinFnGraphCycNodes}\\
\textit{Number of tail nodes }\quad & n-\sqrt{\pi n/2}+1
\tag{iii}\label{T:BinFnGraphTailNodes}\\
\textit{Number of terminal nodes }\quad & n/2
\tag{iv}\label{T:BinFnGraphTerminal}\\
\textit{Number of image nodes }\quad & n/2
\tag{v}\label{T:BinFnGraphImage}
\end{align}
\end{theorem}
In part~(\ref{T:BinFnGraphComp}), $\gamma$ represents the Euler
constant which is approximately 0.57721566.  The highlights of the
proofs as they differ from those in~\cite{RandMap} follow.
\begin{proof}
As
in~\cite{RandMap}, the following bivariate generating functions
need to be defined with parameter $u$ marking the elements of
interest.  The generating functions for the number of components,
number of cyclic nodes and number of terminal nodes are
respectively:
\begin{align}
\xi_1(u,z) &= \exp{\left(u\ln{\frac{1}{1-zb(z)}}\right)}\label{EQ:xiComp}\\
\xi_2(u,z) &= \frac{1}{1-uzb(z)}\label{EQ:xiCYclic}\\
\xi_3(u,z) &= \frac{1}{\sqrt{1-2uz^2}}\label{EQ:xiTerminal}
\end{align}

Imitating the methods
in~\cite{RandMap}, the mean value generating function, $\Xi(z)$, is
found by taking the partial derivative of $\xi(u,z)$ with respect to
$u$ and evaluating at $u=1$.  This yields the following results
\begin{align}
\Xi_1(z) &=
\frac{1}{1-zb(z)}\ln{\left(\frac{1}{1-zb(z)}\right)}\label{EQ:XIComp}\\
\Xi_2(z) &= \frac{zb(z)}{(1-zb(z))^2}\label{EQ:XICyclic}\\
\Xi_3(z) &= \frac{z^2}{(1-2z^2)^{3/2}}.\label{EQ:XITerminal}
\end{align}

The forms in the statement of the theorem follow by expanding around
the singularity $z=1/\sqrt{2}$, applying singularity analysis as
in~\cite{RandMap}, and normalizing parts~(\ref{T:BinFnGraphComp})
and~(\ref{T:BinFnGraphCycNodes}) by~\eqref{EQ:NumGraphs}
and~(\ref{T:BinFnGraphTerminal}) by~\eqref{EQ:NumGraphsExact}.
Parts~(\ref{T:BinFnGraphTailNodes}) and~(\ref{T:BinFnGraphImage})
follow from parts~(\ref{T:BinFnGraphCycNodes})
and~(\ref{T:BinFnGraphTerminal}) respectively since the respective
pairs must sum to $n$.  Also note that 
part~(\ref{T:BinFnGraphTerminal}) can also be derived in an 
elementary fashion from the definition of the binary functional graph.
\end{proof}

The asymptotic values for the average length of cycles and tails as
seen from a random point in the graph are also interesting.  The
asymptotic forms of these values are given in
Theorem~\ref{T:BinFunAvg}.
\begin{theorem}\label{T:BinFunAvg}
The expected values for the cycle size and tail length as seen from
a random node in a random binary functional graph of size $n$ are
asymptotic to
\begin{align}
\textit{Average cycle length }\quad & \sqrt{\pi
n/8}\tag{i}\label{T:BinFunAvgCycle}\\
\textit{Average tail length }\quad & \sqrt{\pi
n/8}\tag{ii}\label{T:BinFunAvgTail}
\end{align}
\end{theorem}
\begin{proof}
In order to calculate the average cycle length and average tail
length, the generating functions must be manipulated to account
for each node in the cycle or tail.  This can be done by using the
same methods as in the previous proof, but on the component
function and taking an additional derivative with respect to $z$
to weight each cycle and tail by the nodes involved.  Multiplying
again by $z$ replaces the factor lost in the differentiation and
by $1/(1-b(z))$ cumulates over all of the components.  This
strategy is used to prove the result for average cycle size
in~\cite{RandMap}.  More background on the method can be found
there.

Marking the appropriate elements, performing a singularity analysis of
the two generating functions and normalizing by
${2^{n/2}}/{(n\sqrt{\pi n/2})},$ as done in the previous theorems,
leads to the statement of the theorem.  The additional factor of $n$
in the denominator is needed to compensate for the fact that the
parameters were estimated across all nodes in the graph and the goal
is to determine them from any single random node in the graph.
\end{proof}

The final parameters that needs to be calculated are the average
maximum cycle length and the average maximum tail length.  

\begin{theorem}\label{T:BinFunMax}
The asymptotic forms for the expected sizes of the largest cycle and 
the largest tail in a random binary functional graph of size
$n$, as $n\rightarrow\infty$, are
\begin{align}
\textit{Largest cycle }\quad &\sqrt{\frac{\pi
n}{2}}\int_0^\infty\left[1-\exp{\left(-\int_v^\infty
e^{-u}\frac{du}{u}\right)}\right]\,dv\approx 0.78248\sqrt{n}
\tag{i} \label{T:BinFunMaxCycle}\\
\textit{Largest tail }\quad & \sqrt{2\pi n} \ln 2 -3 + 2 \ln 2\approx
1.73746\sqrt{n} - 1.61371
\tag{ii}\label{T:BinFunMaxTail}
\end{align}
\end{theorem}
\begin{proof}
The proof for part~\eqref{T:BinFunMaxCycle} result follows
precisely the methods of~\cite{RandMap} with substitution of the
proper generating function $f$, and is therefore omitted.

The proof for part~\eqref{T:BinFunMaxTail} follows a combination
of~\cite[Theorem~6]{RandMap} and~\cite [Sections 3--5]{Average}.  Let
$b^{[h]}(z)$ be the exponential generating function for the number
of binary trees with height at most $h$ and $f^{[h]}(z)$ be the
exponential generating function for the number of binary functional
graphs with maximum tail length less than or equal to $h$, so that (as in
Equations~(41) and~(42) of \cite{RandMap})
$$f^{[h]}(z) = \frac{1}{1-zb^{[h]}(z)}$$
and
$$b^{[h+1]}(z) = z + \frac{1}{2} z \left(b^{[h]}(z)\right)^{2}, \qquad 
b^{[0]}(z) = z.$$

Now, as in~\cite[Proposition~2]{Average}, note that
$$b(z) -
b^{[h+1]} = \frac{1}{2} z \left(b(z) - b^{[h]}(z)\right) \left(b(z) +
b^{[h]}(z)\right)$$
so if we let
$$e_{h}(z) = \frac{b(z) - b^{[h]}(z)}{2b(z)}.$$
then
$$e_{h+1}(z) = (1 - \sqrt{1 - 2z^{2}})e_{h}(z) ( 1 - e_{h}(z)).$$

Now we want to approximate $e_{h}(z)$ with a function of $h$ and some 
$\epsilon(z)$.  If we let 
$\epsilon = \sqrt{1 - 2z^{2}}$
then we have 
$$e_{j+1} = (1 - \epsilon)e_{j} ( 1 - e_{j}); \qquad e_{-1}=2.$$ 
This is essentially the same recursion as 
in~\cite{Average}, and
as in~\cite[Lemma~5]{Average}, we can then 
``normalize'' and ``take inverses'' to get the approximation
\begin{equation} \label{morehope}
    e_{h} \approx 
    \frac{(1-\epsilon)^{h+1}\epsilon}{1-(1-\epsilon)^{h+1}}.
\end{equation}
(The details of the error bounds proceed as in~\cite{Average}; we 
omit them here.)

The generating function associated to the average maximum tail
length is (as in Equation~(43) of \cite{RandMap}) 
$$\Xi(z) = \sum_{h \geq 0} \left[ \frac{1}{1-zb(z)} - 
\frac{1}{1-zb^{[h]}(z)} \right]$$
and we proceed as in Equation~(51) of \cite{RandMap} to write
$$\Xi(z) = \frac{2zb(z)}{1-zb(z)} \sum_{h \geq 0} 
\frac{e_{h}(z) }{1-zb(z)+2e_{h}(z)zb(z)} .$$

Putting this entirely in terms of $\epsilon$ and $h$, and shifting the
index of summation for convenience, we can write
\begin{equation}    
    \label{like52}
 \Xi(z)\approx
\frac{2(1-\epsilon)}{\epsilon} \sum_{h \geq 1} 
\frac{{(1-\epsilon)^{h}}}{1+(1-2\epsilon){(1-\epsilon)^{h}}}.
\end{equation}
We approximate the sum with an integral, using Euler-Maclaurin 
summation.

Taking the integral and noting that $\ln(1-\epsilon) \sim -\epsilon$
as $\epsilon \to 0$, we finally get:
\begin{equation} 
    \label{like55}
 \Xi(z)\approx
\frac{2(1-\epsilon)}{\epsilon^{2} (1-2\epsilon)} \ln (2-3\epsilon + 2 
\epsilon^{2}).
\end{equation}

The next step is to substitute $\epsilon = \sqrt{1 - 2z^{2}}$ 
into~\eqref{like55} and do the singularity analysis, which gives us 
the statement of the theorem.
\end{proof}

\section{Observed Results}\label{S:Observed}
 In~\cite{Holden}, heuristics and observed values for
the number of small cycles (fixed points and two-cycles) in graphs
of the type investigated here are given.  Our methods build on this
to generate experimental data for the parameters described by the
theoretical predictions in Section~\ref{S:Theory}. The method of
data collection was straightforward. A prime was chosen as the
modulus and then for each $g\in\set{1,2,3,...,p-1}$, the
corresponding map or permutation was generated.  The results were
then computed as averages over all $p-1$ graphs observed.  The
permutations and binary functional graphs were noted and their
results were also tabulated separately. In this manner, the data can
be examined in its complete form over all graphs and individually
over the permutations and binary functional graphs.  The generation
and analysis of each of the graphs was handled by C++ code written
by the first author.

The primes chosen for these calculations were \begin{align}100043
&= 2\cdot 50021 + 1,\notag\\
100057 &= 2^3\cdot 3\cdot 11\cdot 379 + 1\text{, and}\notag\\
106261 &= 2^2\cdot 3\cdot 5\cdot 7\cdot 11\cdot 23 + 1.\notag
\end{align}
The total number of graphs, permutations and binary functional
graphs can be computed using Theorem~\ref{T:m-ary} and are shown
in Table~\ref{Table:NumGraphs}.
\begin{table}[h]
\begin{center}
\renewcommand{\arraystretch}{1.25}
\begin{tabular}{|l|c|c|c|}
\hline
& {\bf 100043} & {\bf 100057} & {\bf 106261}\\
\hline
\textbf{Permutations} & 50020 & 30240 & 21120\\
\hline
\textbf{Binary Functional Graphs} & 50020 & 15120 & 10560 \\
\hline
{\bf Total Functional Graphs} & 100042 & 100056 & 106260\\
\hline
\end{tabular}
\caption{The number of permutations, binary functional graphs and
total functional graphs associated with $p=100043$, $p=100057$,
and $p=106260$.}\label{Table:NumGraphs}
\end{center}
\end{table}
The combined results of all functional graphs will be examined
first in Section~\ref{SS:CombinedObserved} where the observed
results will be compared to the theoretical framework for random
functional graphs given in Theorem~\ref{T:FnGraphs}. In
Section~\ref{SS:PermutationObserved}, the observed results for the
permutations will be compared to the theoretical results given in
Theorem~\ref{T:Permutation}.  Finally, the observed results for
the binary functional graphs will be examined in
Section~\ref{SS:BinaryObserved}.  Theorems~\ref{T:BinFunDirect}
through~\ref{T:BinFunMax} will provide the theoretical
predictions for these values.  Since the terminal nodes and tail
nodes can be directly computed from the image nodes and cyclic
nodes, including them in the collected data does not add any
insight.  For this reason, they have both been excluded from the
analysis conducted in the following sections.
Appendix~\ref{A:Extreme} gives some of the interesting extremal
data such as the longest cycle observed for each prime.

\subsection{Combined Results}\label{SS:CombinedObserved}
It would seem that by combining better than one hundred thousand
functional graphs generated by~\eqref{EQ:General}, the results
would tend toward a random functional graph. Theorem~\ref{T:m-ary}
shows that the modular exponentiation function imposes some
structure onto the functional graphs, but especially if $p-1$ has
a complex factorization, the large number of graphs might be thought 
to approach a lack of structure.  However, as
Table~\ref{Table:CombinedResults} clearly shows, these graphs are
not tending toward a random functional graph.
\begin{table}[h]
\begin{center}
\renewcommand{\arraystretch}{1.25}
{\small\begin{tabular}{|l|c|c|c|c|c|c|} \hline
&\multicolumn{2}{|c|}{\bf 100043} & \multicolumn{2}{|c|}{\bf
100057} & \multicolumn{2}{|c|}{\bf 106261}\\
\cline{2-7} &  {\bf
Observed} & {\bf Error} & {\bf Observed} & {\bf Error}& 
{\bf Observed} & {\bf Error}\\\hline 
{\bf Components} & 9.235 & 44.481\% & 7.603 & 18.947\% & 
6.742 & 4.983\%\\\hline 
{\bf Cyclic Nodes} & 50271.600 & 12578.567\% & 30399.400 & 7574.478\% &
21268.600 & 5110.130\%\\\hline 
{\bf Image Nodes} & 75029.000 & 18.644\% & 47838.800 & 24.363\% & 
69435.300 & 3.374\%\\\hline 
{\bf Avg Cycle} & 25088.934 & 12557.883\% & 15249.500 & 7593.148\% &
10629.500 & 5103.529\%\\\hline 
{\bf Avg Tail} & 197.951 & 0.130\% & 114.215 & 42.380\% & 
92.590 & 54.674\%\\ \hline 
{\bf Max Cycle} & 31320.700 & 12555.466\% & 19027.821 & 7587.860\% & 
13259.600 & 5098.564\%\\
\hline
{\bf Max Tail} & 271.408 & 50.613\%& 217.842 & 60.363\%& 
202.581 & 64.232\% \\\hline
\end{tabular}}
\caption{The observed results for the three primes over all functional graphs
generated and the corresponding percent errors.}\label{Table:CombinedResults}
\end{center}
\end{table}

\subsection{Permutation Results}\label{SS:PermutationObserved}

The results in Section~\ref{S:Theory} and Section~\ref{SS:CombinedObserved} imply that the
graphs should be split based on the value of $m$, or the possible
in-degrees of each node.   The results of
looking at only the values of $g$ that were a primitive root modulo
$p$ (permutation graphs) can be found in Table~\ref{Table:Permutation}.
\begin{table}[h]
\begin{center}
\renewcommand{\arraystretch}{1.25}
{\small\begin{tabular}{|l|c|c|c|c|c|c|} \hline
&\multicolumn{2}{|c|}{\bf 100043} & \multicolumn{2}{|c|}{\bf 100057} & \multicolumn{2}{|c|}{\bf 106261}\\\cline{2-7}
 &  {\bf Observed} & {\bf Error} & {\bf Observed} & {\bf Error}& {\bf Observed} & {\bf Error}\\\hline
{\bf Components} & 12.081 & 0.083\% & 12.054 & 0.306\% & 12.126 &  0.205\%\\\hline
{\bf Avg Cycle} & 49980.551 & 0.082\% & 50191.352 & 0.326\% & 53105.104 & 0.048\%\\\hline
{\bf Max Cycle} & 62395.488 & 0.102\% & 62627.745 & 0.256\% & 66245.807 & 0.144\%\\\hline
\end{tabular}}
\caption{The observed results for the three primes over the
permutations and the corresponding percent
errors.}\label{Table:Permutation}
\end{center}
\end{table}

The percent error here
is nearly zero in every instance.  This seems to indicate that there
are no obvious structural differences between a random permutation
and a permutation generated by the process used here.

\subsection{Binary Functional Graph
Results}\label{SS:BinaryObserved}

The binary functional graphs should prove more interesting than
the permutations examined in the previous section.  Unlike
permutations, binary functional graphs do not appear to have been
previously studied in detail.  The statistics derived from the
binary functional graphs and the error when compared to the
results derived in Section~\ref{SS:SafePrime} can be found in
Table~\ref{Table:BinaryResults}.
\begin{table}[h]
\begin{center}
\renewcommand{\arraystretch}{1.25}
{\small\begin{tabular}{|l|c|c|c|c|c|c|} \hline
&\multicolumn{2}{|c|}{\bf 100043} & \multicolumn{2}{|c|}{\bf 100057} & \multicolumn{2}{|c|}{\bf 106261}\\\cline{2-7}
&  {\bf Observed} & {\bf Error} & {\bf Observed} & {\bf Error}& {\bf Observed} & {\bf Error}\\\hline
{\bf Components} & 6.389 & 0.047\% & 6.364 & 0.437\% & 6.370 & 0.810\%\\\hline
{\bf Cyclic Nodes} & 395.303 & 0.029\% & 395.858 & 0.105\% & 408.433 & 0.217\%\\\hline
{\bf Image Nodes} & 50021 & 0\% & 50028 & 0\% & 53130 & 0\%\\\hline
{\bf Avg Cycle} & 198.319 & 0.056\% & 197.766 & 0.230\% & 202.651 & 0.795\%\\\hline
{\bf Avg Tail} & 197.961 & 0.125\% & 197.550 & 0.339\% & 202.422 & 0.907\%\\\hline
{\bf Max Cycle} & 247.261 & 0.094\% & 247.302 & 0.082\% & 256.986 & 0.754\%\\\hline
{\bf Max Tail} & 541.827 & 1.115\% & 549.588 & 1.145\% & 566.370 & 1.744\%\\\hline
\end{tabular}}
\caption{The observed results for the three primes over all binary
functional graphs generated and the corresponding percent
errors.}\label{Table:BinaryResults}
\end{center}
\end{table}

The number of image nodes came out exactly as expected and predicted
by Theorem~\ref{T:image}.  However, in many other cases the results
were nearly as good.  The relative
size of the error follows the number of binary functional graphs for
each prime.  This is especially worth noting for $p=100043$ which
has over fifty thousand binary functional graphs while 100057 and
106261 have approximately fifteen thousand and ten thousand
respectively.  Since having more graphs appears to push the results
closer to those derived in Section~\ref{SS:SafePrime}, this seems to
further support the claim that the results hold for any binary
functional graph produced by our mapping.

\section{Conclusions and Future Work}\label{S:Conclusions}

The transformation used here to generate functional graphs and
permutations is an exceedingly important transformation in
cryptography.  If the output of the function were to fall into a
predictable pattern, it could be an exploitable flaw in many
algorithms considered secure today.  For instance, the average cycle
length seems particularly important for pseudorandom bit generators since, in many cases,
it relates directly to the predictability of the pseudorandom bit generator. As
Theorem~\ref{T:m-ary} demonstrates, the use of~\eqref{EQ:General}
repeatedly forces a non-trivial structure onto the graphs generated.
This is certainly worth investigating as any imposed structure may
be open to an exploit.

The advantage of using a safe prime is that every non-trivial
graph can be analyzed by the theoretical framework laid out in
this paper.  Their use is also very prevalent in cryptographic
applications.  As mentioned above, the pseudorandom bit generator
specified in~\cite{PRBG} requires the use of a safe prime to
defend against other attacks.  However, the methods used for
binary functional graphs in Section~\ref{SS:SafePrime} can and
should be extended to larger values of $m$.  In an ideal case,
they should be extended in the general case for an $m$-ary graph
that can be specified by
\begin{center}
\begin{tabular}{ll}
FunctionalGraph & = set(Components)\\
Component & = cycle(Node*Set(Tree, cardinality = $m-1$))\\
Tree & = Node + Node*set(Tree, cardinality = $m$)\\
Node & = Atomic Unit
\end{tabular}
\end{center}
The associated generating functions for these functional graphs
would be
\begin{align*}
f(z) &= e^{c(z)}\\
c(z) &= \ln\left(1-\frac{z}{(m-1)!}t^{m-1}(z)\right)^{-1}\\
t(z) &= z + \frac{z}{m!}t^m(z)
\end{align*}
where $f(z)$ is the exponential generating function associated to
the functional graphs, $c(z)$ is the exponential generating
function associated to the connected components and $t(z)$ is
associated to the trees.  The methods in
Section~\ref{SS:SafePrime} could also be extended to obtain values
for additional parameters such as the maximum tail length.

This paper has focused on the graphs generated when the modulus is
prime.  In practice, though, this is not always the case.  For
this reason, it could be worthwhile to attempt to extend the type
of analysis done here to a composite modulus.

While the data generated for this project appears to confirm that the
graphs do tend toward the shape and structure of a random graph of the
appropriate type, no data was collected on the distribution of the
different parameters.  This data could help to give a clearer picture
of how closely individual graphs may be expected to exhibit the
characteristics of a random graph, especially given the observation
that primes with a larger number of binary functional graphs seem to 
conform better to prediction on the average.  The methods used 
in~\cite{DistribHeights} would seem to be potentially helpful here.

\appendix

\section{Extremal Data}\label{A:Extreme}
For $p=100043$, the longest cycle observed was 100042 which
occurred for two different values of $g$.  They were $g=20812$ and
$g=94034$.  The longest tail had a length of 1448 and was observed
when $g=89339$.  There were five instances where the graphs
contained no cycles longer than one which occurred for $g=1$,
72116, 91980, 95997, and 100042.

The graphs generated by $p=100057$ had an overall longest cycle of
100052 when $g=58303$.  The longest tail observed was 1589 when
$g=18115$.  There were also 26 different values of $g$ that
produced a graph that did not have a cycle longer than one.

The largest cycle observed in graphs generated using $p=106261$
was 106257 when $g=102141$.  The longest tail was 35822 when
$g=1480$.  There were 92 different values of $g$ that produced
graphs with no cycles longer than a fixed point.

\end{document}